\documentclass[12pt]{article}  % list options between brackets
\usepackage{times,amsmath,amsfonts,amsthm,amssymb,eucal}
\usepackage{array}
\usepackage{mathrsfs}
\usepackage{amsmath}
\usepackage{setspace}
\usepackage{fullpage}
\usepackage{fancyhdr}
\usepackage{setspace}
\usepackage{amssymb}
\usepackage{epsfig}
\usepackage{graphicx}
\usepackage{makeidx}
\usepackage{color}
\usepackage{multirow}
\usepackage{epsfig}
\usepackage{graphicx}
\usepackage{makeidx}
\usepackage{gensymb}
\usepackage{enumerate}
\usepackage{multirow}
 \usepackage{setspace}
\usepackage{amssymb}
 \usepackage{setspace}

\usepackage[english]{babel}
\usepackage{amsmath,amsthm}
\usepackage{amsfonts}

\usepackage{multirow}

\usepackage{subfigure}
\usepackage{epsfig}

\usepackage{lscape}
\usepackage{geometry}

\newtheorem{theorem}{Theorem}
\newtheorem{proposition}{Proposition}
\newtheorem{lemma}{Lemma}
\newtheorem{corollary}{Corollary}
\newtheorem{definition}{Definition}
\def\Theorem{\begin{theorem}\sl}
\def\EndTheorem{\end{theorem}}
\def\Proposition{\begin{proposition}\sl}
\def\EndProposition{\end{proposition}}
\def\Lemma{\begin{lemma}\sl}
\def\EndLemma{\end{lemma}}
\def\Corollary{\begin{corollary}\sl}
\def\EndCorollary{\end{corollary}}
\def\Definition{\begin{definition}\sl}
\def\EndDefinition{\end{definition}}
\theoremstyle{remark}

\def\mref#1{(\ref{#1})}

\numberwithin{equation}{section}

\begin{document}
%%%%%%%%%%%%%%%%%%%%%%%%%%%%%%%%%%%%%%%%%%%%%%%%%%%%%%%%%%%%%%%%%%%%%%%%%%%%%%%%%%%%%%%%%%%%%%%%%%%%%%%%%%%%%%%%%%%%%%%%%%%%
\title{ \textbf{Two-sample Bayesian nonparametric goodness-of-fit test}}   % type title between braces

\author{Luai Al Labadi \thanks{{\em Address for correspondence}: Luai Al Labadi, Department of Statistical Sciences, University of Toronto, Toronto, Ontario M5S 3G3, Canada. E-mail: luai.allabadi@utoronto.ca},  Emad  Masuadi \thanks{Emad  Masuadi, Research Unit, Department of Medical Education, College of Medicine, King Saud Bin Abdulaziz University for Health Sciences, Riyadh, KSA. E-mail: masuadie@ksau-hs.edu.sa}\& Mahmoud Zarepour\thanks{M. Zarepour, Department of Mathematics and Statistics,
University of Ottawa, Ottawa, Ontario, K1N 6N5, Canada. E-mail: zarepour@uottawa.ca} }
\date{\today}    % type date between braces
\maketitle

\pagestyle {myheadings} \markboth {} {Two-sample Bayesian nonparametric GOF test}

\begin{abstract}
%Testing the difference between two data samples is of a particular interest in statistics.  Precisely, given two samples $\mathbf{X}=X_1,\ldots,X_{m_1} \overset {i.i.d.} \sim F$  and $\mathbf{Y}=Y_1,\ldots,Y_{m_2} \overset {i.i.d.} \sim G$, with $F$ and $G$ being unknown continuous cumulative distribution functions, we wish to test the null hypothesis $\mathcal{H}_0:~F=G$. In this paper, we propose an effective and  convenient Bayesian nonparametric  approach to assess the equality of two unknown distributions.  The method is based on the Kolmogorov distance and approximate samples from the Dirichlet process centered at  the standard normal distribution and a concentration parameter  1. Our results show that the proposed test is robust with respect to any prior specification of the Dirichlet process. We provide  simulated  examples to illustrate the workings of the method. Overall, the proposed method performs perfectly in many cases.
In recent years, Bayesian nonparametric statistics has gathered extraordinary attention. Nonetheless, a relatively little amount of work has been expended on Bayesian nonparametric hypothesis testing. In this paper, a novel Bayesian nonparametric approach to the two-sample problem is established. Precisely, given two samples $\mathbf{X}=X_1,\ldots,X_{m_1}$ $\overset {i.i.d.} \sim F$  and $\mathbf{Y}=Y_1,\ldots,Y_{m_2} \overset {i.i.d.} \sim G$, with $F$ and $G$ being unknown continuous cumulative distribution functions, we wish to test the null hypothesis $\mathcal{H}_0:~F=G$.   The method is based on the Kolmogorov distance and approximate samples from the Dirichlet process centered at  the standard normal distribution and a concentration parameter  1. It is demonstrated that the proposed test is robust with respect to any prior specification of the Dirichlet process. A power comparison with several well-known tests is  incorporated. In particular, the proposed test dominates the standard Kolmogorov-Smirnov test in all the cases examined in the  paper.


\vspace{9pt}

%\vspace{9pt} \noindent { \textbf{Mathematics Subject Classification (2000):}} 62F15; 62G10; 60B10.
\end{abstract}

%%%%%%%%%%%%%%%%%%%%%%%%%%%%%%%%%%%%%%%%%%%%%%%%%%%%%%%%%%%%%%%%%%%%%%%%%%%%%%%%%%%%%%%%%%%%%%%%%%%%%%%%%%%%%%%%%%%%%%%%%%%%

\section{Introduction}
\label{intro}
Two-sample comparison is a common problem in statistics. Namely, given two samples $\mathbf{X}=X_1,\ldots,X_{m_1} \overset {i.i.d.} \sim F$ and $\mathbf{Y}=Y_1,\ldots,Y_{m_2} \overset {i.i.d.} \sim G$, with $F$ and $G$ being unknown continuous cumulative distribution functions, the problem  is to decide whether $F=G$. For instance, in medical studies, one may want to assess the efficiency of a new drug in two groups of patients. See Borgwardt and Ghahramani (2009) for more interesting examples in different disciplines.

The objective of this paper is to describe a Bayesian nonparametric procedure for the above situation. Our method is based on approximate samples from the Dirichlet process (Ferguson, 1973) with the standard normal base measure  and a concentration parameter of unity. Next, the Kolmogorov distance is used to examine if the two distributions are equal or not.

Bayesian nonparametrics is a fast developing area in statistics. Nevertheless, there has been relatively little amount of work has been expended on Bayesian nonparametric hypothesis testing. Most of the work includes goodness-of-fit tests for one-sample problems. Two standard nonparametric Bayesian approaches for one-sample goodness-of-fit tests can be found  in the literature. The first approach consists of embedding the proposed model in the null hypothesis into a larger family of models  (the alternative family). Following this step, a prior is placed on the alternative family. Then, the Bayes factor  of the null hypothesis to the alternative is computed. For example, Carota and Parmigiani (1996), and Florens, Richard, and Rolin (1996) used a Dirichlet process prior for the alternative distribution. McVinish, Rousseau, and Mengersen (2009) considered  mixtures of triangular distributions. Another form of the prior, the P\'olya tree process (Lavine, 1992), was suggested by Berger and Guglielmi (2001). The  second approach for one-sample goodness-of-fit tests is based on placing a prior on the true distribution generating the data. For this test,  the distance between the posterior distribution and the proposed one is measured.  Muliere and Tardella (1998), Swartz (1999), Al Labadi and Zarepour (2013a, 2014b) considered the Dirichlet process and applied the Kolmogorov distance to test continuous distributions.  Viele (2000) used the Dirichlet process and the Kullback-Leibler distance to test only discrete distributions. Explicit expressions for calculating the different types of distance between the Dirichlet process and its base measure were derived in Al Labadi and Zarepour (2014b). On the other hand, Hsieh (2011) used the P\'olya tree prior and the Kullback-Leibler distance to test continuous distributions.  As for two-sample tests, Holmes, Caron, Griffin,  and Stephens (2015)  developed a way to compute the Bayes factor for testing the null hypothesis through the marginal likelihood of the data with P\'olya tree priors centered either subjectively or using an empirical procedure. Under the null hypothesis, they modeled the two samples to  come from a single random measure distributed as a P\'olya tree, whereas under the alternative hypothesis the two samples come from two separate P\'olya tree random measures. Ma and Wong (2011) allowed the two distributions to be generated jointly through optional coupling of a P\'olya tree prior. Borgwardt and Ghahramani (2009) discussed two-sample tests based on Dirichlet process mixture models and derived a formula to compute the Bayes factor in this case. Generalizations of the Bayes factor approach based on P\'olya tree priors to censored and multivariate data were proposed  by Chen and Hanson (2014). Note that, the two-sample Bayesian nonparametric tests based on the distance approach are not found in the literature. Thus, the method proposed in this paper is considered the first endeavor in this direction.

%In the non-Bayesian literature, nonparametric hypothesis testing for two sample difference has a long history and rich literature, as well. Several procedures have been proposed. Well-known examples are Komlogorov-Smirnov test, Cram\'er-von Mises test, and Wilcoxon test. See Lehmann(i.e., the concentration parameter and the base measure) and Recently, this problem has also been investigated from a Bayesian nonparametric perspective using the P\'olya tree process (Holomes, Caron, Griffin,  and Stephens, 2009).

This paper is structured as follows. In Section 2, we recall the definition of the Dirichlet process and some of its relevant properties.  In Section 3, we describe our method to test the equality of two unknown distributions. In Section 4, illustrative examples and simulation results are included. In Section 5, we empirically compare the power of the proposed test to several well-known tests. Some properties of the proposed approach  are discussed in Section 6. Finally, some concluding remarks are made in Section 7.

\section{The Dirichlet Process} \label{intro}
In this  section, we   introduce some preliminary information about the Dirichlet process. The Dirichlet process, formally introduced in Ferguson (1973), is the most well-known and widely used prior in Bayesian nonparametric inference. Consider a space $\mathfrak{X}$ with a $\sigma-$algebra $\mathcal{A}$ of subsets of $\mathfrak{X}$. Let $H$ be a fixed probability measure on $(\mathfrak{X},\mathcal{A})$ and $a$ be a positive number. Following Ferguson (1973), a random probability
measure $P=\left\{P(A)\right\}_{A \in \mathcal{A}}$ is called a Dirichlet process on $(\mathfrak{X},\mathcal{A})$ with parameters
$a$ and $H$, if for any finite measurable partition $\{A_1, \ldots, A_k\}$ of $\mathfrak{X}$, the joint distribution of the vector
$\left(P(A_1), \ldots\,P(A_{k})\right)$ has the Dirichlet distribution with parameters $(a H(A_1), \ldots,$ $aH(A_k)),$ where $k\ge 2$. We assume that if $H(A_j)=0$, then $P(A_j)=0$ with a probability one. If $P$ is a Dirichlet process with parameters $a$ and $H,$ we write $P\sim {DP}(a, H).$ The parameter $a$ is known as the \emph{concentration parameter} and the probability measure $H$ is called the \emph{base (centering) measure} of $P.$

An attractive feature of the Dirichlet process is the conjugacy property. If  $X_{1},\ldots, X_{m}$ is a  sample from $P\sim DP(a, H)$, then the posterior distribution of $P$ given $ X_{1},\ldots, X_{m}$ coincides with the distribution of the Dirichlet process
 with parameters $a^{*}$ and $H^{*}$, where
\begin{equation}
 a^{*}=a+m \quad \text{and} \quad H^{*}_{m}=\frac{a }{a+m}H+\frac{m}{a+m}\frac{\sum_{i=1}^{m}\delta_{{X}_i}}{m}.\label{eq2.2}
\end{equation}
Here and throughout this paper, $\delta_X$ denotes the Dirac measure at $X$, i.e. $\delta_X(A)=1$ if $X \in A$ and $0$ otherwise. We also use a ``$*$" as a superscript to denote posterior quantities. James, Lijoi, and Pr\"unster (2006) showed that the Dirichlet process is  the only normalized random measure with independent increments that enjoys the conjugacy property. Notice that, the posterior base distribution $H^{*}$ is a convex combination of the base distribution  and the empirical distribution. The weight associated with the prior base distribution $H$ is
 proportional to $a$, while the weight
 associated with the empirical distribution is proportional to the number of observations $m$.
The posterior base distribution $H^*$  approaches the prior base measure $H$ for large values of $a.$ On the
 other hand, for small values of $a $, $H^*$ is close to the empirical distribution. The consistency property of the Dirichlet process has been studied in detail in Goshal (2010). Similar to the  frequentist's empirical process, as $m\to \infty$, Lo (1987) showed that  the centered and scaled Dirichlet process $\sqrt{m}\left({P}^*_m-H^*_m\right)$
converges to a Brownian bridge on $D[0,1]$ with respect to the Skorohod topology. Lo (1987) applied his result to establish asymptotic validity of the Bayesian bootstrap. See also James (2008) and Al Labadi and  Zarepour (2013b).  The distributional functionals of the Dirichlet process appear,  for instance, in  Cifarelli and  Regazzini (1990), Regazzini, Guglielmi, and Di Nunno  (2002), Lijoi and Regazzini (2004) and James (2005, 2006).

Ferguson (1973) proposed a series representation  as an alternative definition for the Dirichlet process. Also see Ferguson and Klass  (1972). Specifically, let $(E_k)_{k \ge 1}$ be a sequence of i.i.d. random variables  with an exponential distribution of mean 1 and $\Gamma_i=E_1+\cdots+E_i$.
Let $(\theta_i)_{i \geq 1}$  be a sequence of i.i.d. random variables with common distribution $H$, independent of $(\Gamma_i)_{i\ge 1}$.
Define
\begin{equation}
P(\cdot)=\sum_{i=1}^{\infty} {\frac{L^{-1}(\Gamma_i)}{\sum_{i=1}^{\infty}{{L^{-1}(\Gamma_i)}}}\delta_{\theta_i}(\cdot)},\label{eq2.4}
\end{equation}
where $L(x)=a\int_{x}^{\infty}t^{-1}e^{-t}dt, x>0,$ and
$L^{-1}(y)=\inf\{x>0: L(x) \geq y\}.$ Then  the random probability measure $P$ is a Dirichlet process with parameters $a$ and $H$.

From (\ref{eq2.4}), it follows clearly that a realization of the Dirichlet process is
a discrete probability measure. This is true even when the base measure is absolutely continuous. This fact was noted by Ferguson (1973), and Blackwell and MacQueen (1973). Note that, although the Dirichlet process is  discrete with probability one, this discreteness is no more troublesome than the discreteness of the empirical process. By imposing the weak topology, the support for the Dirichlet process is quite large. Specifically, the support for the Dirichlet process is the set of all probability measures whose support is contained in the support of the base measure. This means  if the support of the base measure is $\mathfrak{X}$, then the space of all probability measures is the support  of the Dirichlet process. For example, if we have  a normally distributed base measure, then the Dirichlet process can choose any probability measure. See Ferguson (1973) and Ghosh and Ramamoorthi (2003) for further discussion about the support of the Dirichlet process. In practice,\ it is difficult to work with~\mref{eq2.4}  because there is no tractable form  for the L\'evy measure $L$ and  determining the random weights in the sum requires the computation of  an infinite sum. Recently, Zarepour and Al Labadi (2012) derived an efficient  approximation of the Dirichlet process with monotonically decreasing weights. Specifically, let $X_n$  be a random variable with a $\text{Gamma}(a/n,1)$ distribution. Define
\begin{equation}
G_n(x)=\Pr(X_n>x)={\frac{1}{\Gamma(a/n)}\int_{x}^{\infty}e^{-t}t^{a/n-1}dt} \nonumber
\end{equation}
and
\begin{equation}
G^{-1}_n(y)=\inf\left\{x:G_n(x)\ge y\right\}. \nonumber
\end{equation}
Let $(\theta_i)_{1 \leq i \leq n}$ be a sequence of i.i.d. random variables
with values in $\mathfrak{X}$ and common distribution $H$, independent
of $(\Gamma_i)_{1 \leq i \leq {n+1}}.$ Then, as $n\to \infty$,
\begin{equation}
P_n=\sum_{i=1}^{n} \frac
{{G_n^{-1}\left(\frac{\Gamma_i}{\Gamma_{n+1}}\right)}}{\sum_{i=1}^{n}{G_n^{-1}\left(\frac{\Gamma_i}{\Gamma_{n+1}}\right)}}\delta_{\theta_i}
\label{eq11}
\end{equation}
converges almost surely to $P$, defined by (\ref{eq2.4}). Zarepour and Al Labadi (2012) and Al Labadi and Zarepour (2014a) demonstrated that the convergence rate
of the representation (\ref{eq11}) is empirically faster than several existing representations, including Bondesson (1982), Sethuraman (1994), and Wolpert and Ickstadt (1998).

Based on representation (\ref{eq11}), the following algorithm outlines the steps required to generate a sample from the approximate Dirichlet process with parameters $a$ and $H$:

\noindent \textbf{Algorithm A: Simulating an approximation for the Dirichlet process process}
\begin{enumerate} [(1)]
\item Fix a relatively large positive integer $n$.
\item Generate  $\theta_i\overset{\text{i.i.d.}}\sim H$ for $i=1,\ldots,n.$
\item For $i=1,\ldots,n+1,$ generate $E_i$ from an exponential distribution with  mean 1, independent of $\left(\theta_i\right)_{1\le i \le n}$ and let $\Gamma_i=E_1+\cdots+E_i.$
\item For $i=1,\ldots,n,$ compute $G_n^{-1}\left({\Gamma_i}/{\Gamma_{n+1}}\right),$  which is  simply the quantile function of the $\text{Gamma}(\alpha/n,1)$ distribution evaluated at  $1-{\Gamma_i}/{\Gamma_{n+1}}.$
\item  Use representation (\ref{eq11}) to find an approximate sample of the Dirichlet process.
 \end{enumerate}

\section{A Bayesian Nonparametric  goodness-of-fit Test}
In this section, we consider the two-sample problem described in the Introduction, where two i.i.d. samples are observed and the problem is to test if  the two underlying distributions are different. Specifically, given two samples $\mathbf{X}=X_1,\ldots,X_{m_1} \overset {i.i.d.} \sim F$ and $\mathbf{Y}=Y_1,\ldots,Y_{m_2} \overset {i.i.d.} \sim G$ with $F$ and $G$ being unknown continuous cumulative distribution functions, we want to test the null hypothesis $\mathcal{H}_0:~F=G$. The approach is based on measuring the Kolmogorov distance  between the posterior distribution of the Dirichlet process given the first sample and the posterior distribution of the Dirichlet process given the second sample.  Next, we compare whether the  distance  is small or large. Since the posterior distribution of the Dirichlet process converges uniformly to the actual distribution generating the data as the sample size gets large (Ferguson, 1973;  Goshal, 2010; Al Labadi and Zarepour, 2013b), $\mathcal{H}_0$ is rejected whenever the distance is large (see also Lemma \ref{BSP3} in Section 6 for more discussion). On the other hand, $\mathcal{H}_0$ is not rejected if the distance is small. Two issues must be considered: (1) how to select the  parameters for the Dirichlet process and (2) how to conclude whether the resulting distance  is large or small. As  for the first issue, we choose the base measure to be the standard normal distribution and  the concentration parameter to be 1. We show in Section 6 of the current paper that the proposed test is robust  with respect to the prior specification of the Dirichlet process (i.e., the concentration parameter and the base measure). To address the second issue, we introduce first the Kolmogorov distance. Let $P_{n_1}$ and $Q_{n_2}$ be two discrete distributions with corresponding jump points $(U_k)_{1 \leq k \leq n_1}$ and $(V_k)_{1 \leq k \leq n_2}$. The Kolmogorov distance between $P_{n_1}$ and $Q_{n_2}$, denoted by $d(P_{n_1},Q_{n_2}):=d$, is
\begin{equation}
d(P_{n_1},Q_{n_2})=\sup_{x\in \mathbb{R}}\left|P_{n_1}\left((-\infty,x]\right)-Q_{n_2}\left((-\infty,x]\right)\right|:=\sup_{x\in \mathbb{R}}|P_{n_1}(x)-Q_{n_2}(x)|,\nonumber
\end{equation}
where, here and throughout the paper, we use the same notation for the probability measure and  its corresponding cumulative distribution function. The above distance can be simplified (for programming convenience) to
\begin{equation}
d(P_{n_1},Q_{n_2})=\max_{1\le i \le n_1+n_2 }|P_{n_1}(Z_i)-Q_{n_2}(Z_i)|,\label{d800}
\end{equation}
where $(Z_k)_{1 \leq k \leq {n_1+n_2}}$ are the combined jump points. In symbols, $Z_k=U_k, k=1,\ldots,n_1$ and $Z_{n_1+k}=V_k, k=1,\ldots,n_2$. Thus, for each $k$, $k=1,\ldots,n_1+n_2$, we compute $|P_{n_1}(Z_i)-Q_{n_2}(Z_i)|$  and set $d$ to be the largest of these values.

In our approach,  we  set $P_{n_1}=P^{*}_{n_1,m_1}$ and $Q_{n_2}=Q^{*}_{n_2,m_2}$ in (\ref{d800}), where $P^{*}_{n_1,m_1}$ is  an approximation of the posterior distribution of the Dirichlet process   given the first sample and $Q^{*}_{n_2,m_2}$ is an approximation of the posterior distribution of the Dirichlet process  given the second sample. Small values of $d$ indicate  evidence in favor of $\mathcal{H}_0$. To determine whether  $d$  is large or small, the (prior) Kolmogorov distance between two prior  distributions of the Dirichlet process is computed. Henceforth, $d_0$ denotes the Komogorov distance between two prior distributions. We take the base measure for each  prior distribution to be the standard normal distribution, where the concentration parameter of the first prior is $1+m_1$ and for the second prior it is  $1+m_2$. This setup of the prior distributions of the  Dirichlet process guarantees that  any change between the prior distance and posterior distance is only due to the difference between the two samples. Then we calculate the $95\%$ prediction interval by  deleting the lowest and highest $2.5\%$ of the values of $d_0$. We set $U$ to be  the upper bound of the $95\%$ prediction interval of $d_0$. It follows that, we reject (do not reject) $\mathcal{H}_0$ if the mean of the values of $d$ is greater (smaller)  than $U$.  It is straightforward to construct tables for values of $U$ with different sample sizes and  significance levels.  For convenience, for sample sizes less than or equal to 20, values of $U$ are reported in Table 10 in the Appendix. On the other hand, for sample sizes greater than or equal to 20, values of $U$ can be approximated by the following formula:
\begin{equation}
U\approx1.41\sqrt{ \frac {1} {m_1}+ \frac {1}{m_2}} \label{reg1}.
\end{equation}
Formula (\ref{reg1}) is derived via a regression  of $U$ from the simulation of different sample sizes, where the coefficient of determination $R^2$ for the regression equation is more than $99.9\%$. Note that the value of $U$ is determined based on the distance between two Dirichlet processes. Such distance   depends only on the cumulative weights of  each Dirichlet process. That is,  $H$ is not playing any role in determining the distance. Thus, the value of $U$ does not depend on $H$. It follows that Table 10 and formula (\ref{reg1}) are valid for any $H$. See also Section 6 for more discussion about this point.

The following algorithm summarize the steps required for a Bayesian nonparametric  goodness-of-fit test for two samples:
\vspace{3 mm}

\noindent \textbf{Algorithm B:  Bayesian nonparametric test for two samples}
\begin{enumerate}
\item [(1)] Set the base measure $H$ of the Dirichlet process to be the standard normal distribution and the concentration parameter $a$  to  1.
\item [(2)] Use Algorithm A to generate a random sample from an approximation of the posterior Dirichlet process $P^{*}_{n_1,m_1}$, given the first sample. Here $m_1$ is the  size of the first sample.
\item [(3)] Use Algorithm A to generate a random sample from an approximation of the posterior Dirichlet process $Q^{*}_{n_2,m_2}$, given the second sample. Here $m_2$ is the  size of the second sample.
\item [(4)] Compute $d\left(P^{*}_{n_1,m_1},Q^{*}_{n_2,m_2}\right)$, as  defined  in (\ref{d800}).
\item [(5)] Repeat steps (2)-(4) to obtain $r$ i.i.d. samples of $d\left(P^{*}_{n_1,m_1},Q^{*}_{n_2,m_2}\right)$. For large $n_1$, $n_2$ and $r$, the empirical distribution of these values is an approximation to the distribution of $d(P^{*}_{m_1},Q^{*}_{m_2})$, where $P^{*}_{m_1}$ is the posterior distribution of the Dirichlet process given the first sample and $Q^{*}_{m_2}$ is the posterior distribution of the Dirichlet process given the second sample.
  \item [(6)] Calculate $U$, the upper bound of $95\%$ prediction interval confidence interval as follows: (alternatively,  use either Table 10 or  formula (\ref{reg1}))
  \begin{enumerate}
 \item [(i)] Repeat the above steps (1)-(5) to calculate the distance $d_0$ between prior distributions of the Dirichlet process. The base measure is the standard normal distribution, while the concentration parameters  for the first and the second samples are $1+m_1$ and $1+m_2$, respectively.
 \item [(ii)] Sort the values of $d_0$ and set $U$ to be the maximum value after deleting $2.5\%$ of the largest values of $d_0$.
\end{enumerate}
  \item [(7)] If the mean of the distance $d$ is less than $U$, then there is a sufficient evidence not to reject $\mathcal{H}_0$. Otherwise, we reject the null hypothesis $\mathcal{H}_0.$
\end{enumerate}

\section{Examples}
In this section,  the work of the proposed method is illustrated through the following examples.

\noindent {\textbf{Example 1.}}  \label{example1}
Consider  samples generated from the following distributions, where each sample is of size 100. These distributions are also considered in Holmes, Caron, Griffin, and Stephens  (2015).
\begin{enumerate}
  \item $\mathbf{X}\sim N(0,1)$ and $\mathbf{Y}\sim N(0,1)$
  \item $\mathbf{X}\sim N(0,1)$ and $\mathbf{Y}\sim N(1,1)$
  \item $\mathbf{X}\sim N(0,1)$ and $\mathbf{Y}\sim N(0,2)$
  \item $\mathbf{X}\sim N(0,1)$ and $\mathbf{Y}\sim 0.5N(-2,1)+0.5N(2,1)$
  \item $\mathbf{X}\sim N(0,1)$ and $\mathbf{Y}\sim t_{3}$
  \item $\mathbf{X}\sim N(0,1)$ and $\mathbf{Y}\sim t_{0.5}$
  \item $\log \mathbf{X}\sim N(0,1)$ and  $\log \mathbf{Y}\sim N(1,1)$
  \item  $\log \mathbf{X}\sim N(0,1)$ and  $\log \mathbf{Y}\sim N(0,2)$,
  \end{enumerate}
where $N(\mu,\sigma)$ is the normal distribution with mean $\mu$ and standard deviation $\sigma$ and $t_r$ is the $t$ distribution with $r$ degrees of freedom.
In Algorithm B, we set $H=N(0,1)$,  $a=1$, $n_1=n_2=1000$, $r=2000$, and  $m_1=m_2=100$. The results are reported in Table 1. Thus, we reject the null hypothesis whenever the mean of $d$ is greater than $U=0.2$. We also compare our results with standard (frequentist) goodness-of-fit tests  such as the Kolmogorov-Smirnov (K-S) test and the Wilcoxon  (Mann-Whitney U) test. To  calculate these tests we have used the codes ``ks.test" and ``wilcox.test" available in R. It follows from Table 1 that the new test performs very well for all cases. Since the Wilcoxon  test assumes that one of the samples must be a shifted version of the other, it is not used for samples 3, 4, 5, 6, and 8. Therefore, using the Wilcoxon  test is not always reasonable in practice.

\begin{table}[htbp]
  \centering
  \caption{Example 1: Bayesian nonparametric test against (frequentist) Kolmogorov-Smirnov test and  Wilcoxon test.}
    \begin{tabular}{lccc}
    \hline
Samples & $d$: Bayesian &  p-value: K-S & p-value: Wilcoxon  \\
    \hline
1	&	0.15	&		0.8127  	&		0.6002 \\
2	&   0.42 	&		0.0000	&	   0.0000\\
3	&	0.26	&	   0.0101  &		-\\
4	&	0.44	&		0.0000  &	   -\\
5	&	0.19	&		0.2106 	&		- \\
6	&	0.31	&		0.0008 	&		-\\
7	&	0.42	&	0.0000	&		0.0000\\
8	&	0.26	&		0.0101	&		-\\
    \hline
$U$	&	0.20	&			&   \\
    \hline
    \end{tabular}%
  \label{tab:addlabel}%
\end{table}%

Figures 1, 2, 3, and 4 provide plot of  5 sample paths for each of the posterior Dirichlet  process given the first sample  and the posterior Dirichlet  process given the second sample. Conclusions similar to that given above can also be drawn from the figures.
\begin{figure}[h!]
\centering
\vspace{-1 cm}
\subfigure[$\mathbf{X}\sim N(0,1)$ and $\mathbf{Y}\sim N(0,1)$]{\epsfig{figure=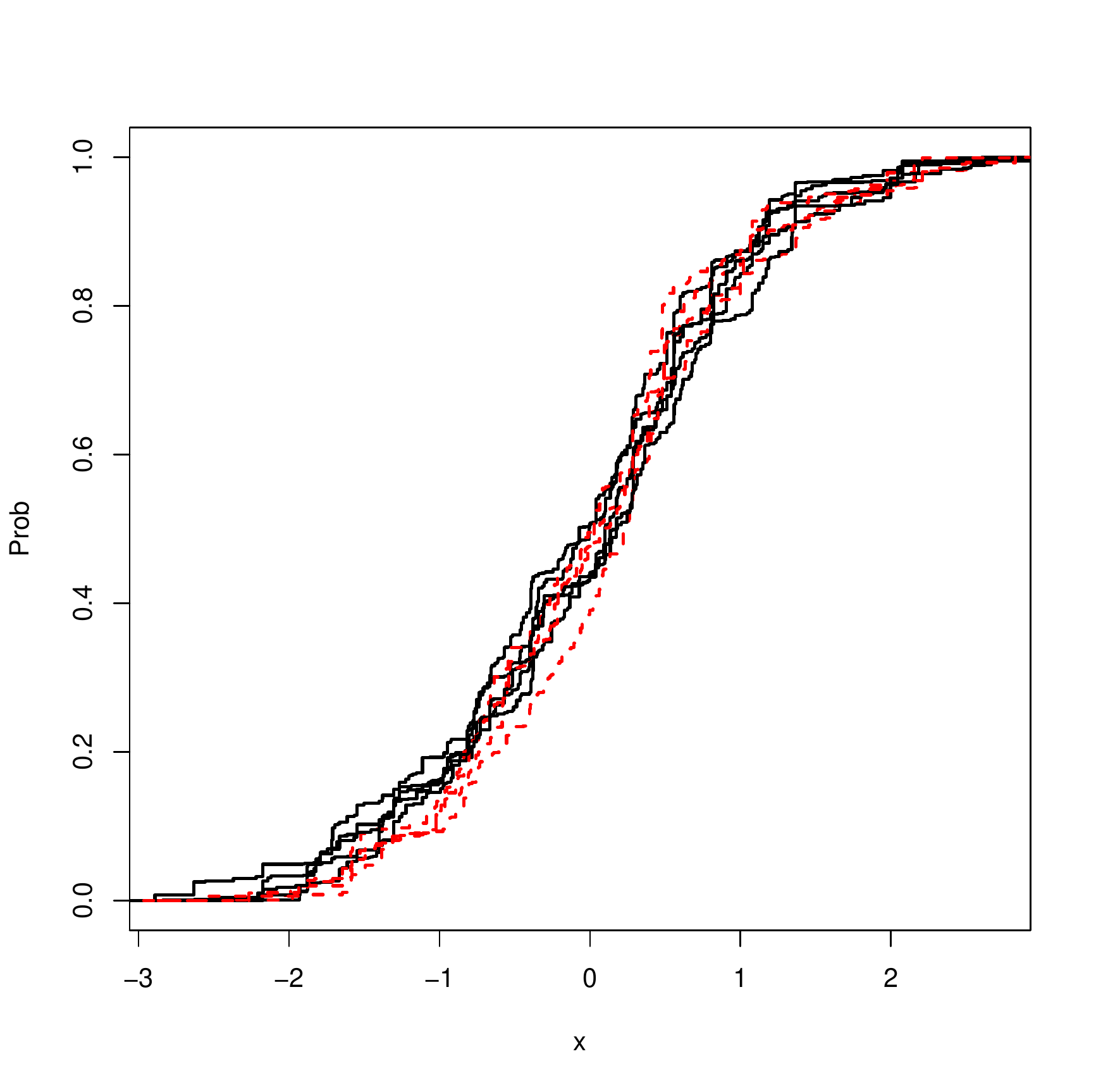,width=4in}}

 \subfigure[$\mathbf{X}\sim N(0,1)$ and $\mathbf{Y}\sim N(1,1)$]{\epsfig{figure=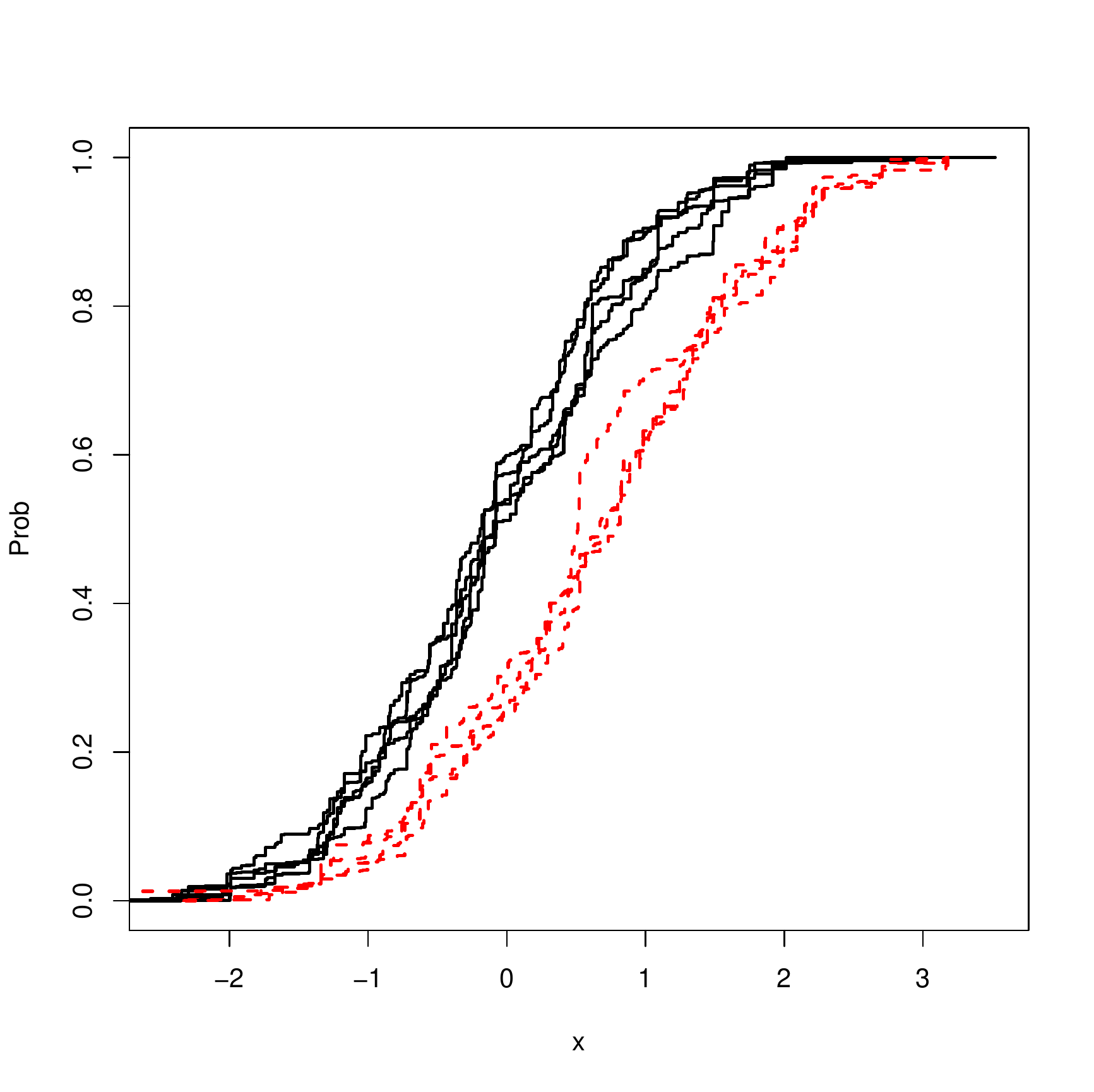,width=4in}}

 \caption{The solid lines represent sample paths of the posterior Dirichlet  process given the first sample and the dashed lines represent sample paths of the posterior Dirichlet process given the second sample.}
 \label{fig:SubF1}
\end{figure}

\begin{figure}[h!]
\centering
\vspace{-1cm}
\subfigure[$\mathbf{X}\sim N(0,1)$ and $\mathbf{Y}\sim N(0,2)$]{\epsfig{figure=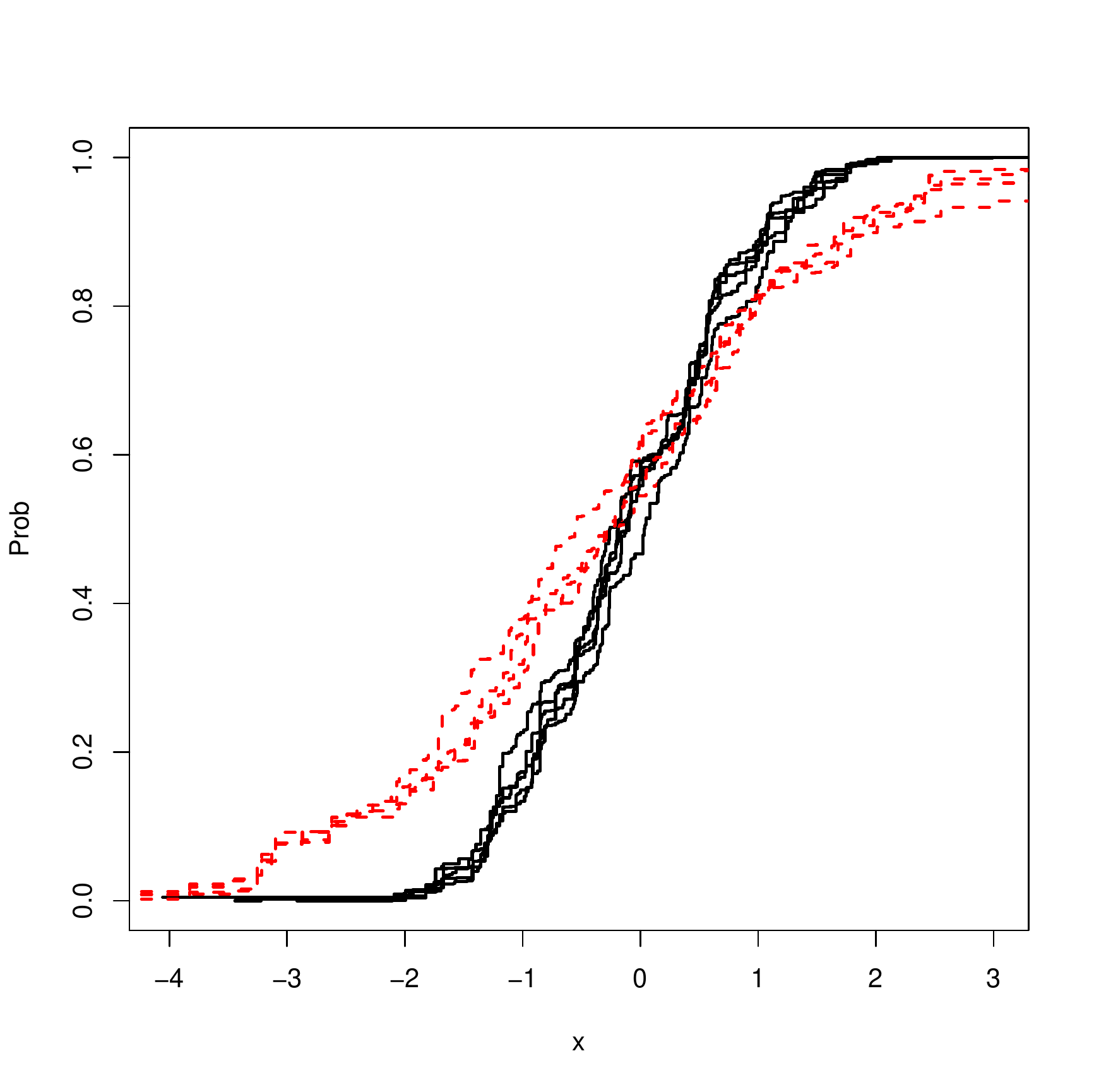,width=4in}}

 \subfigure[$\mathbf{X}\sim N(0,1)$ and $\mathbf{Y}\sim 0.5N(-2,1)+0.5N(2,1)$]{\epsfig{figure=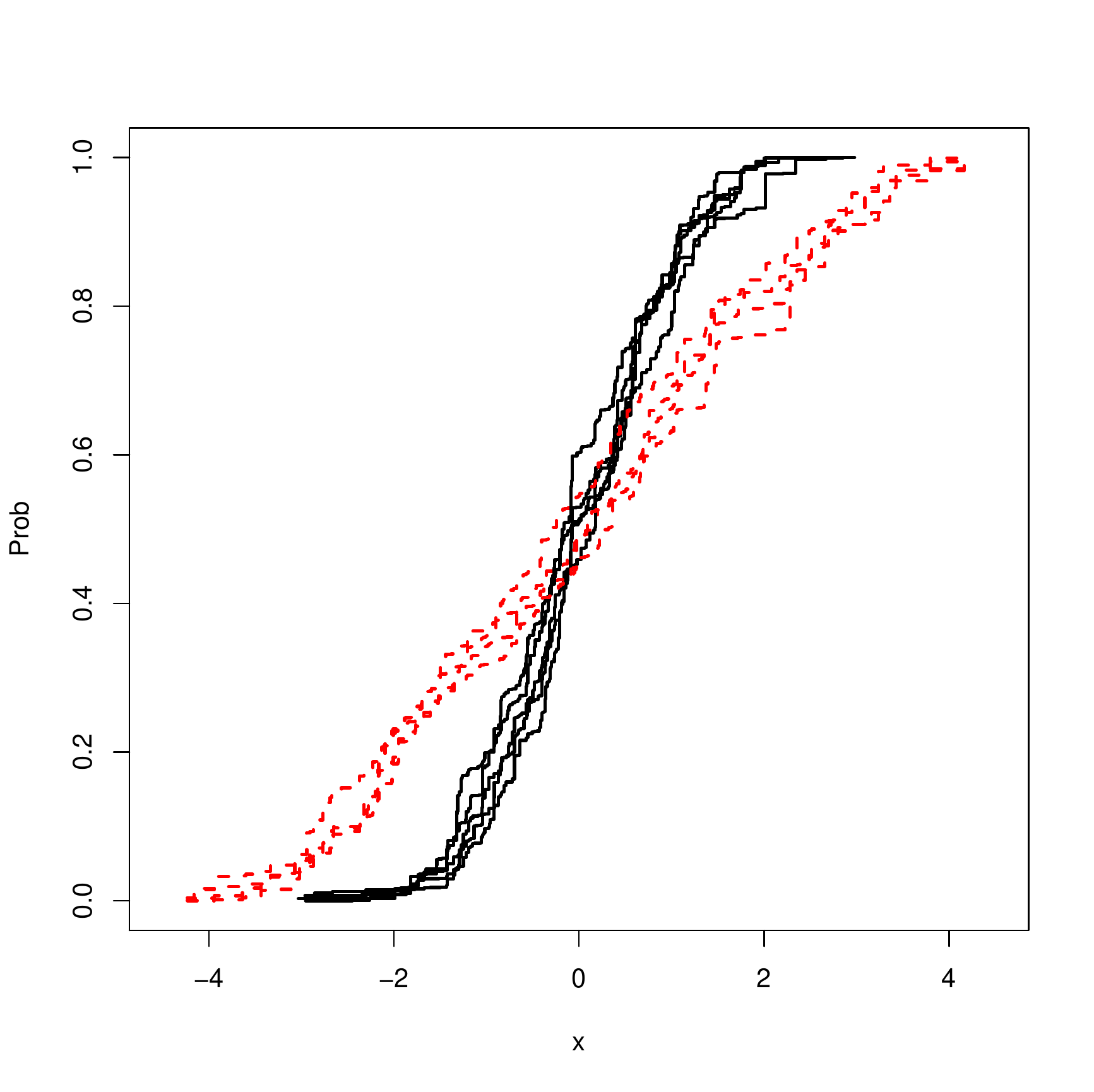,width=4in}}
  \caption{The solid lines represent sample paths of the posterior Dirichlet  process given the first sample and the dashed lines represent sample paths of the posterior Dirichlet process given the second sample.}
 \label{fig:SubF1}
\end{figure}

\begin{figure}[h!]
\centering
\vspace{-1cm}
\subfigure[$\mathbf{X}\sim N(0,1)$ and $\mathbf{Y}\sim t_{3}$]{\epsfig{figure=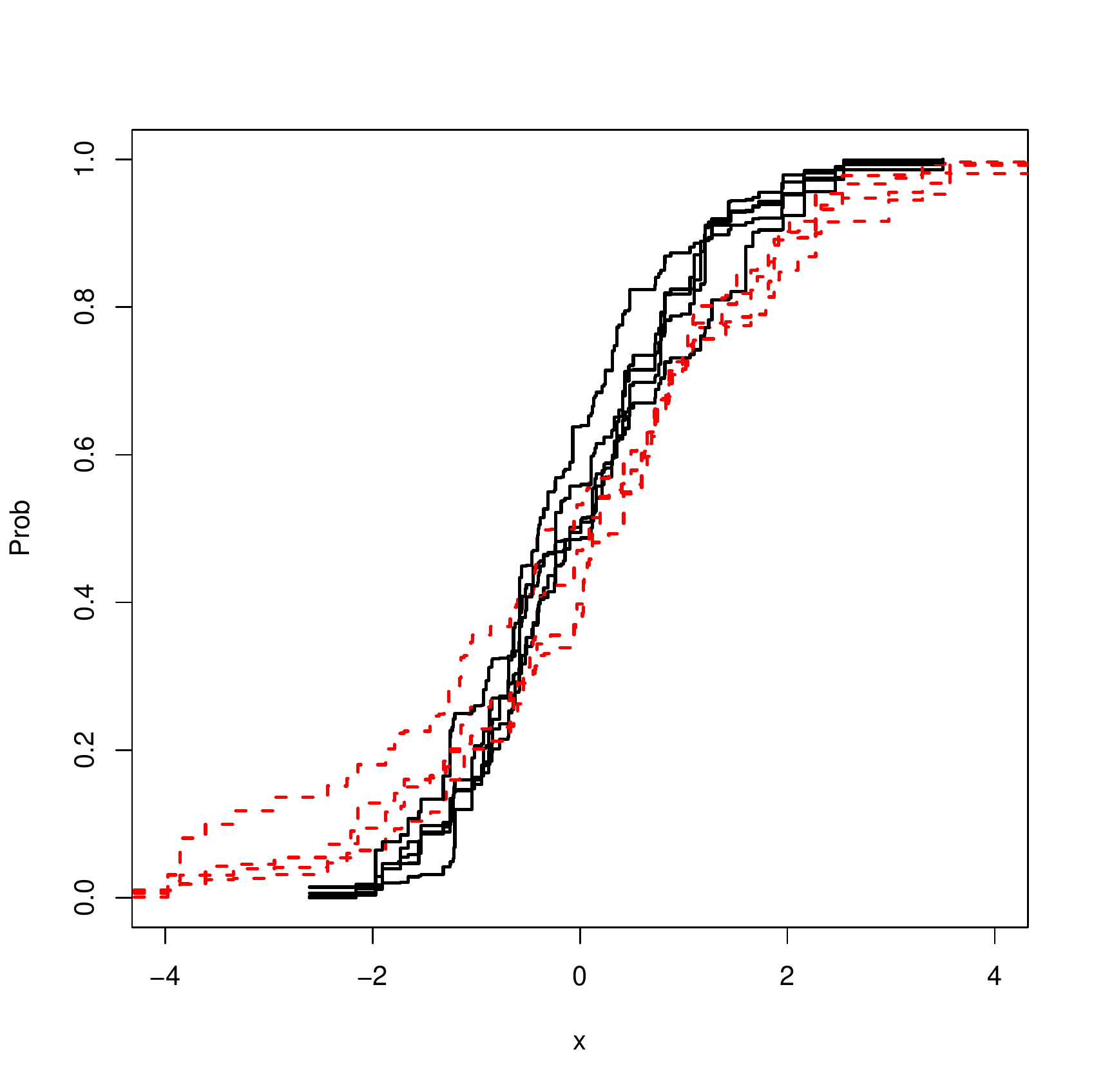,width=4in}}

 \subfigure[$\mathbf{X}\sim N(0,1)$ and $\mathbf{Y}\sim t_{0.5}$]{\epsfig{figure=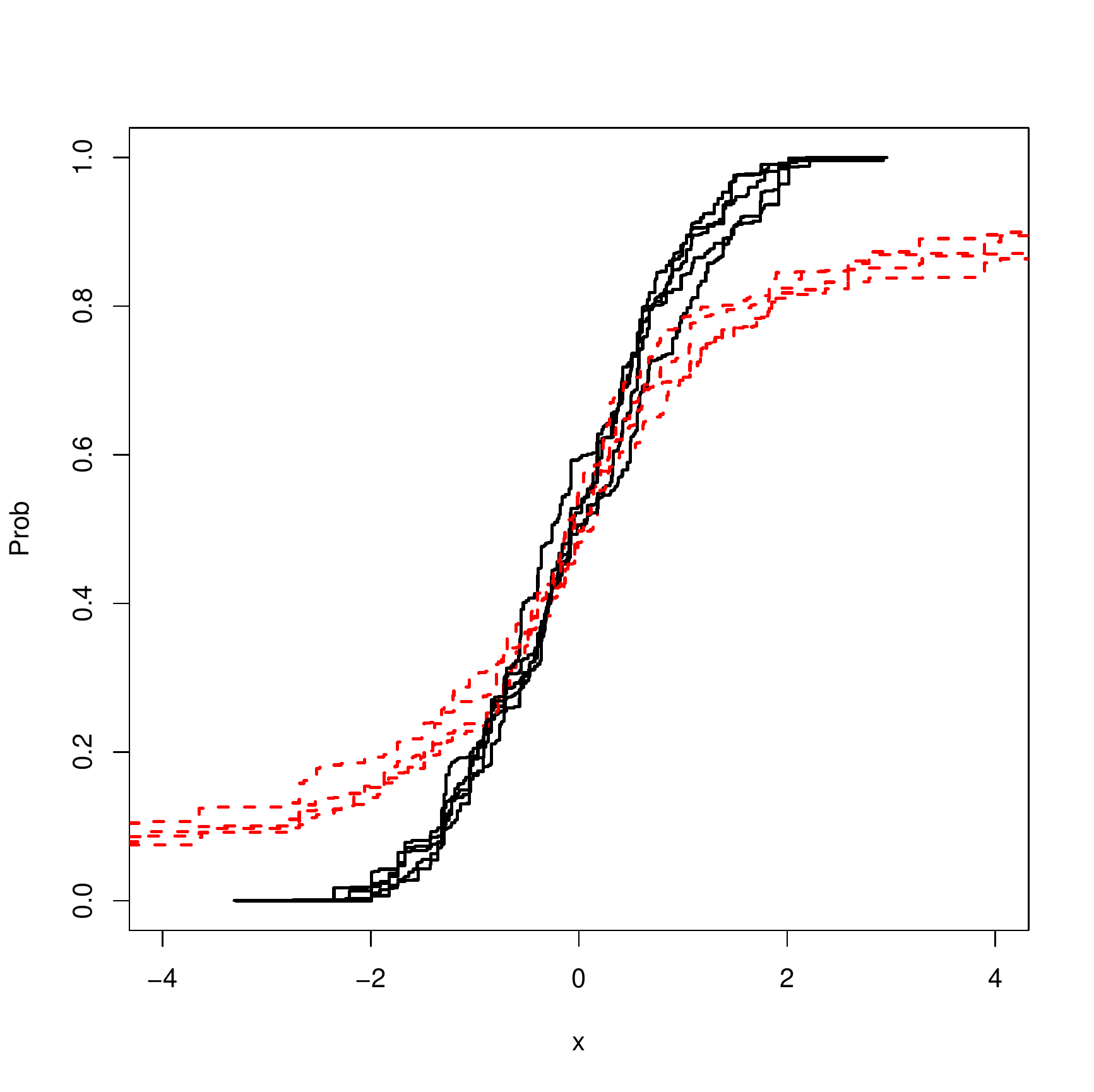,width=4in}}
  \caption{The solid lines represent sample paths of the posterior Dirichlet  process given the first sample and the dashed lines represent sample paths of the posterior Dirichlet process given the second sample.}
 \label{fig:SubF2}
\end{figure}

\begin{figure}[h!]
\centering
\vspace{-1cm}
\subfigure[ $\log \mathbf{X}\sim N(0,1)$ and  $\log \mathbf{Y}\sim N(1,1)$]{\epsfig{figure=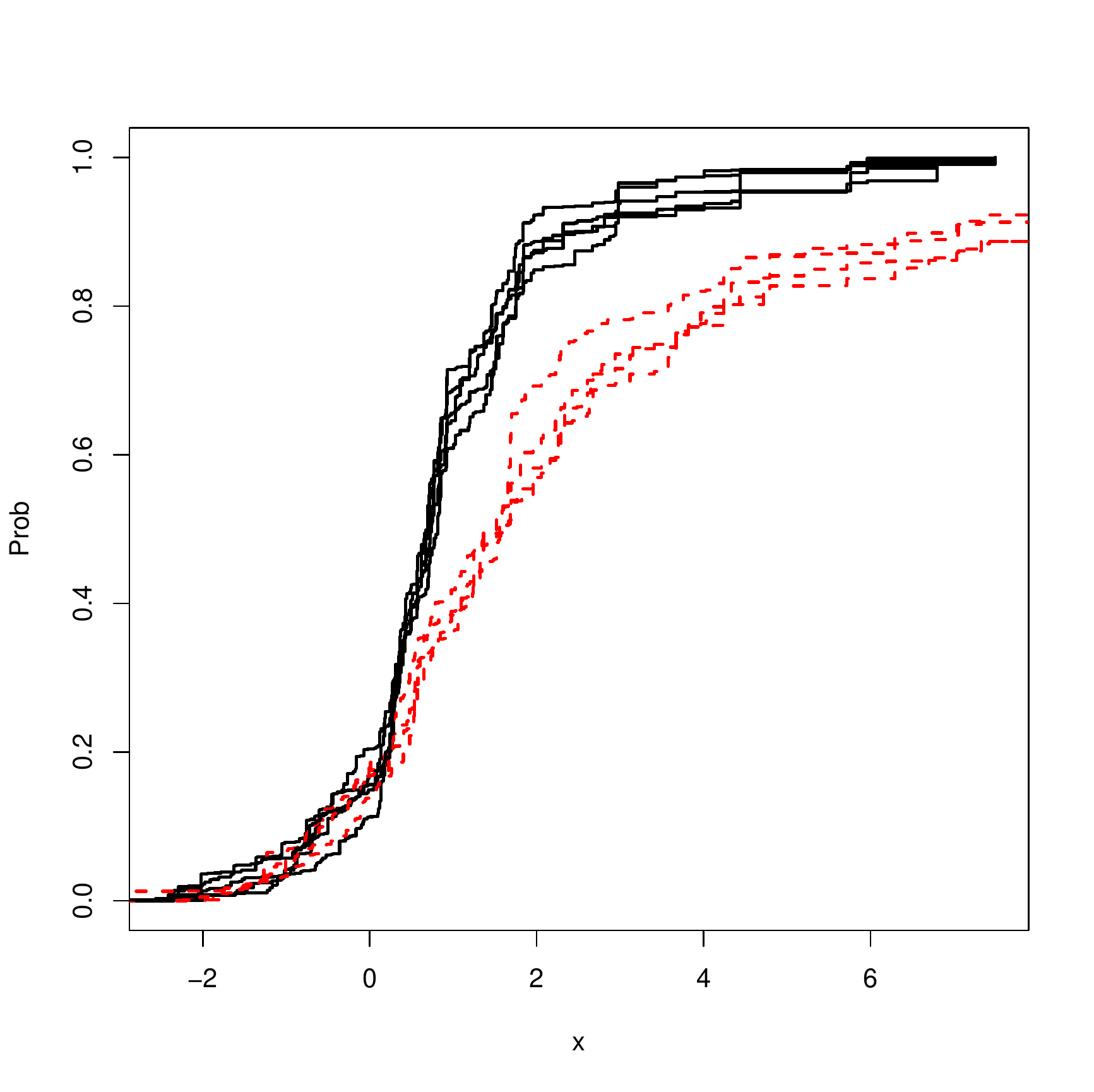,width=4in}}

 \subfigure[$\log\mathbf{X}\sim N(0,1)$ and  $\log \mathbf{Y}\sim N(0,2)$]{\epsfig{figure=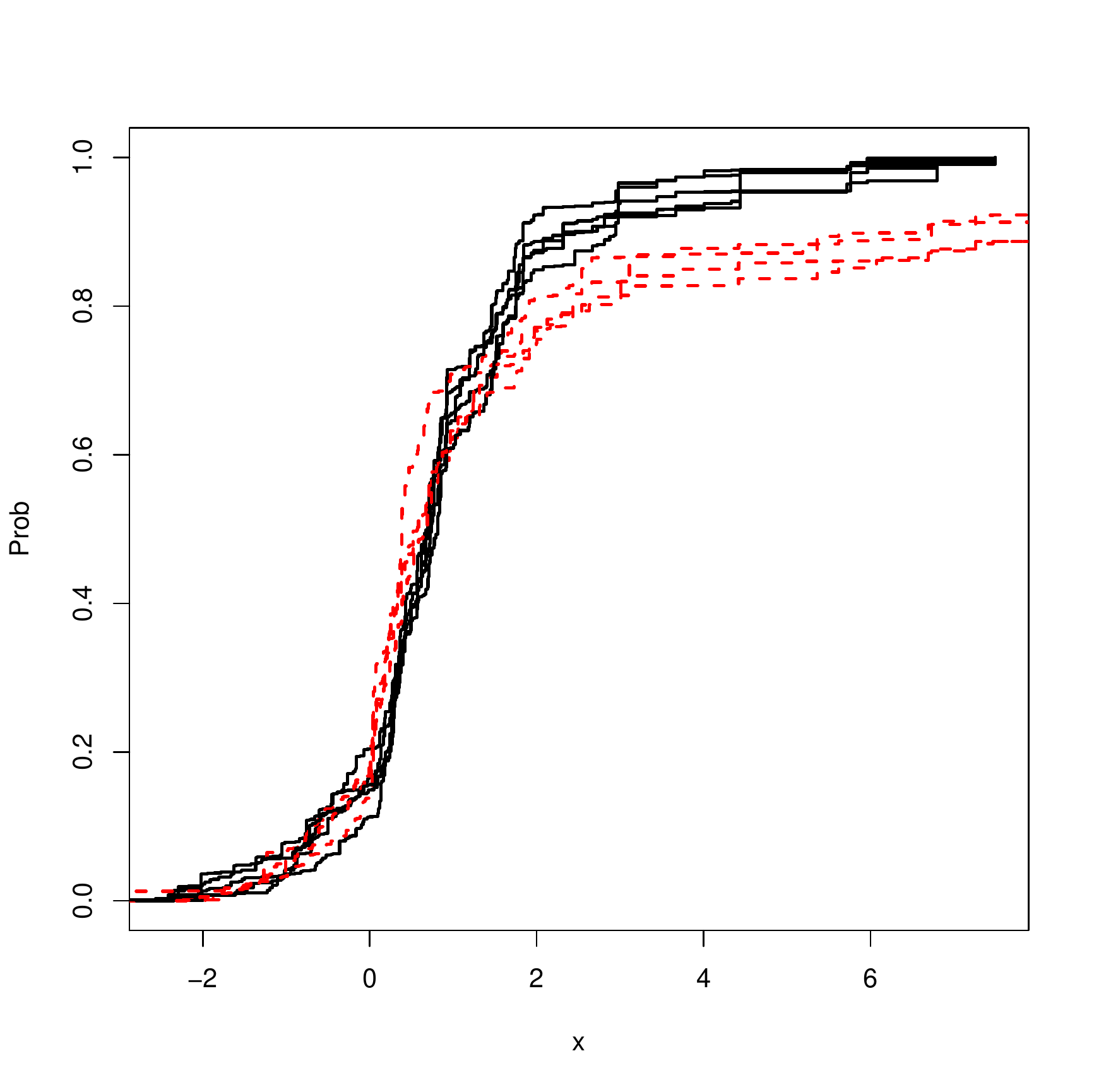,width=4in}}
  \caption{The solid lines represent sample paths of the posterior Dirichlet  process given the first sample and the dashed lines represent sample paths of the posterior Dirichlet process given the second sample.}
 \label{fig:SubF2}
\end{figure}

\noindent {\textbf{Example 2.}}  \label{example2} In this example, we study the performance of the proposed test as the sample size increases. We consider samples from the distributions given in Example 1, cases 1 and 2. The results are summarized in Table 2 and Table 3.

\begin{table}[htbp]
  \centering
  \caption{Example 2: $\mathbf{X}\sim N(0,1)$ and $\mathbf{Y}\sim N(0,1)$.}
    \begin{tabular}{lcccc}
    \hline

Samples & $d$: Bayesian & U& p-value: K-S  & p-value: Wilcoxon  \\
    \hline
$m_1=m_2=5$	    &0.50	    	    &0.69	 	&0.8730          &0.4206\\
$m_1=m_2=10$	&0.41 				&0.56		&0.7869	     	 &0.7959\\
$m_1=m_2=15$	&0.33 				&0.48		&0.9383 		 &0.8381\\
$m_1=m_2=20$	& 0.31 				&0.42		&0.5713 	 	 &0.5683\\
$m_1=m_2=30$	&0.23				&0.36		&0.3929 		 &0.2301\\
$m_1=m_2=50$	&0.14 				&0.28		&0.7166 	 	 &0.5192\\
$m_1=m_2=100$	&0.15				&0.20		&0.5806	    	 &0.3958\\
$m_1=m_2=200$	&0.12 				&0.16		&0.5441 	     &0.7598\\
    \hline
    \end{tabular}%
  \label{tab:addlabel}%
\end{table}%

\begin{table}[htbp]
  \centering
  \caption{Example 2: $\mathbf{X}\sim N(0,1)$ and $\mathbf{Y}\sim N(1,1)$}
\begin{tabular}{lcccc}
    \hline

Samples & $d$: Bayesian & U& p-value: K-S  & p-value: Wilcoxon  \\
    \hline
$m_1=m_2=5$&0.72    	&0.70&0.0794       	 &0.0317\\
$m_1=m_2=10$&0.73   	&0.56&0.0123  	  	 &0.0005\\
$m_1=m_2=15$&0.61		&0.47&0.0077 	 	 &0.0010\\
$m_1=m_2=20$&0.53		&0.42&0.00397 	      &0.0132\\
$m_1=m_2=30$&0.52		&0.36&0.0000			 &0.0000\\
$m_1=m_2=50$&0.50		&0.28&0.0000			 &0.0000\\
$m_1=m_2=100$&0.46		&0.20&0.0000			 &0.0000\\
$m_1=m_2=200$&0.39		&0.16&0.0000			 &0.0000\\
    \hline
    \end{tabular}%
  \label{tab:addlabel}%
\end{table}%

It follows from Table 2  that, in all cases, the null hypothesis is not rejected. On the other hand, from  Table 3, the null hypothesis is  rejected. In both scenarios, the results are consistent with that obtained by the  frequentist tests. Thus, the proposed test works even with sample sizes as small as 5. The power of the proposed test is examined in the next section.

\section{Power Comparison} In this section, we empirically compare the power of the proposed test with the (standard) Kolmogorov-Smirnov test and the Wilcoxon test. We consider the following cases:
\begin{enumerate}
  \item $\mathbf{X}\sim N(0,1)$ and $\mathbf{Y}\sim N(1,1)$
  \item $\mathbf{X}\sim N(0,1)$ and $\mathbf{Y}\sim N(0,2)$
  \item $\mathbf{X}\sim N(0,1)$ and $\mathbf{Y}\sim 0.5N(-2,1)+0.5N(2,1)$
  \item $\mathbf{X}\sim N(0,1)$ and $\mathbf{Y}\sim t_{0.5}$
  %\item  $\log \mathbf{X}\sim N(0,1)$ and  $\log \mathbf{Y}\sim N(0,2)$,
   \item  $\mathbf{X}\sim \text{Exponential}(1)$ and  $ \mathbf{Y}\sim \text{Exponential}(2)$.  [$F_X(x)=1-e^{-x}$ and $F_Y(y)=1-e^{-2y}$]
  \end{enumerate}
The power of each test is estimated  by calculating the proportion of rejecting the null hypothesis out of 1000 samples when the two sampling distributions are different (i.e. $F\neq G$). We consider several sample sizes. The simulation results are summarized in Table 4 to Table 8. The results in the tables indicate that the power of the proposed test is always greater than the power of the  Kolmogorov-Smirnov test. On the other hand, as pointed out in Section 4, Wilcoxon  test assumes that one of the samples must be a shifted version of the other. In practice, such assumption is not known a priori. Thus, Wilcoxon  test is only useful for Case 1 and Case 5. Even in these two  cases, the proposed test is as efficient as the  Wilcoxon  test. See Table 4 and Table 8.

\begin{table}[htbp] %1
  \centering
  \caption{Power comparison from $\mathbf{X}\sim N(0,1)$ and $\mathbf{Y}\sim N(1,1)$.}
    \begin{tabular}{lcccc}
    \hline

Samples & Bayesian & K-S &  Wilcoxon  \\
    \hline
$m_1=m_2=5$	    &0.343	    	    	 	&0.075          &0.235\\
$m_1=m_2=10$	&0.480 						&0.218	     	 &0.484\\
$m_1=m_2=15$	&0.670 						&0.533 		 &0.711\\
$m_1=m_2=20$	& 0.840 						&0.724 	 	 &0.855\\
$m_1=m_2=30$	&0.948						&0.885 		 &0.963\\
$m_1=m_2=50$	&0.995 						&0.982 	 	 &0.995\\
$m_1=m_2=100$	&1						&1	    	 &1\\
$m_1=m_2=200$	&1 						&1 	     &1\\
    \hline
    \end{tabular}%
  \label{tab:addlabel}%
\end{table}%

\begin{table}[htbp] %2
  \centering
  \caption{Power comparison from $\mathbf{X}\sim N(0,1)$ and $\mathbf{Y}\sim N(0,2)$.}
    \begin{tabular}{lcccc}
    \hline

Samples & Bayesian & K-S  &  Wilcoxon  \\
    \hline
$m_1=m_2=5$	    &0.123	    	    	 	&0.01          &0.049\\
$m_1=m_2=10$	&0.108 						&0.022	     	 &0.049\\
$m_1=m_2=15$	&0.125						&0.067		 &0.048\\
$m_1=m_2=20$	&0.217 						&0.124 	 	 &0.054\\
$m_1=m_2=30$	&0.332					& 0.184		 &0.059\\
$m_1=m_2=50$	&0.667					& 0.371		 &0.051\\
$m_1=m_2=100$	&0.979 					&0.800	 	 &0.060\\
$m_1=m_2=200$	&1						&0.997	    	 &0.059\\
    \hline
    \end{tabular}%
  \label{tab:addlabel}%
\end{table}%

\begin{table}[htbp] %3
  \centering
  \caption{Power comparison from $\mathbf{X}\sim N(0,1)$ and $\mathbf{Y}\sim 0.5N(-2,1)+0.5N(2,1)$.}
    \begin{tabular}{lcccc}
    \hline

Samples & Bayesian & K-S  &  Wilcoxon  \\
    \hline
$m_1=m_2=5$	    &0.176	    	    	 	&0.025          &0.056\\
$m_1=m_2=10$	&0.303 						&0.075	     	 &0.065\\
$m_1=m_2=15$	&0.545 						&0.277 		 &0.079\\
$m_1=m_2=20$	& 0.784 						&0.472	 	 &0.074\\
$m_1=m_2=30$	&0.949						&0.762		 &0.075\\
$m_1=m_2=50$	&1 						&0.982	 	 &0.076\\
$m_1=m_2=100$	&1						&1	    	 &0.066 \\
$m_1=m_2=200$	&1 						& 1	     &0.075\\
    \hline
    \end{tabular}%
  \label{tab:addlabel}%
\end{table}%

\begin{table}[htbp] %4
  \centering
  \caption{Power comparison from $\mathbf{X}\sim N(0,1)$ and $\mathbf{Y}\sim t_{0.5}$.}
    \begin{tabular}{lcccc}
    \hline

Samples & Bayesian & K-S  &  Wilcoxon  \\
    \hline
$m_1=m_2=5$	    &0.109	    	    	 	&0.015        &0.046\\
$m_1=m_2=10$	&0.107						&0.028     	 &0.057\\
$m_1=m_2=15$	&0.135 						&0.061 		 &0.057\\
$m_1=m_2=20$	& 0.226						&0.110 	 	 &0.058\\
$m_1=m_2=30$	&0.338						&0.187 		 &0.059\\
$m_1=m_2=50$	&0.762						&0.491 	 	 &0.049\\
$m_1=m_2=100$	&1						&0.943	    	 &0.061\\
$m_1=m_2=200$	&1 						&1 	     &0.051\\
    \hline
    \end{tabular}%
  \label{tab:addlabel}%
\end{table}%

\begin{table}[htbp] %5
  \centering
  \caption{Power comparison from $\mathbf{X}\sim \text{Exponential} (1)$ and  $ \mathbf{Y}\sim \text{Exponential} (2)$.}
    \begin{tabular}{lcccc}
    \hline

Samples & Bayesian & K-S  &  Wilcoxon  \\
    \hline
$m_1=m_2=5$	    &0.173	    	    	 	&0.026          &0.100\\
$m_1=m_2=10$	&0.202 						&0.069     	 &0.209\\
$m_1=m_2=15$	&0.315						&0.226 		 &0.337\\
$m_1=m_2=20$	& 0.446					&0.321	 	 &0.444\\
$m_1=m_2=30$	&0.589						&0.484		 &0.620\\
$m_1=m_2=50$	&0.847						&0.746 	 	 &0.832\\
$m_1=m_2=100$	&0.994 						&0.961 	 	 &0.991\\
$m_1=m_2=200$	&1					&1 	     &1\\
    \hline
    \end{tabular}%
  \label{tab:addlabel}%
\end{table}%

\section{Additional Properties }

In this section, we discuss some additional properties of the proposed method for the goodness-of-fit test. The next lemma shows that, under the null hypothesis, the Kolmogorov distance between posterior distributions given the data converges to zero as sample sizes get large.

\Lemma
\label{BSP3}
If $P^{*}_{m_1}$  is the posterior distribution of the Dirichlet process  given the first and $Q^{*}_{m_2}$ is the posterior distribution of the Dirichlet process  given the second sample. Then, under $\mathcal{H}_0$, we have
\begin{equation}
d(P^{*}_{m_1},Q^{*}_{m_2}) \to 0, \nonumber
\end{equation}
as $m_1,m_2\to \infty$. Recall, $m_i$ represents the sample size of the sample $i$, $i=1,2$.
\EndLemma

\proof It follows from the triangle inequality that
%\begin{equation}
%d(P^{*}_{m_1},Q^{*}_{m_2}) \le d(P^{*}_{m_1},H^{*}_{m_1})+d(Q^{*}_{m_2},H^{*}_{m_2})+d(H^{*}_{m_1},H^{*}_{m_1}), \label{triangle}
%\end{equation}
\begin{eqnarray}
\nonumber d\left(P^{*}_{m_1},Q^{*}_{m_2}\right) &\le& d\left(P^{*}_{m_1},H^{*}_{m_1}\right)+d\left(Q^{*}_{m_2},H^{*}_{m_1}\right)\\ \label{triangle} &\le&  d\left(P^{*}_{m_1},H^{*}_{m_1}\right)+d\left(Q^{*}_{m_2},H^{*}_{m_2}\right)+d\left(H^{*}_{m_1},H^{*}_{m_2}\right),
\end{eqnarray}
where $H^{*}_{m_i}$ is defined in (\ref{eq2.2}) for any continuous base measure $H$, $i=1,2$. The proof of the lemma is complete since, under $\mathcal{H}_0$,  the right hand side of the inequality (\ref{triangle}) converges to zero as $m_1,m_2\to \infty$ (Ferguson, 1973;  Goshal, 2010;   Al Labadi and Zarepour, 2013b).
\endproof

%Let $Z_{(1)}<\ldots <Z_{{(n_1+n_2)}}$ be the order statistics
%of $(Z_k)_{1 \leq k \leq n_1+n_2}$.  Recall that $(Z_k)_{1 \leq k \leq {n_1+n_2}}$ are the combined jump points. That is, $Z_k=U_k, k=1,\ldots,n_1$ and $Z_{n_1+k}=V_k, k=1,\ldots,n_2$. Define $U_{(0)}=V_{(0)}=-\infty$. The next proposition gives an approximation of the Kolmogorov prior distance.
%
%\Lemma
%\label{BSP4} Let $P_{n_1}=\sum_{i=k}^{n_1} p_k\delta_{U_k}$  and $Q_{n_2}=\sum_{i=k}^{n_2} q_k\delta_{V_k}$ with jump
%points $(U_k)_{1 \leq k \leq n_1}$ and $(V_k)_{1 \leq k \leq n_2}$ come from a continuous distribution $H$. We have
%\begin{equation}
%     d(P_{n_1},Q_{n_2})=\max\left\{d^{(1)},d^{(2)},d^{(3)},d^{(4)}\right\},\nonumber
%\end{equation}
%where
%\begin{eqnarray*}
%d^{(1)}=\max_{1\le i \le n_1+n_2}\left|\sum_{k=1}^{i-1} p'_k-\sum_{k=1}^{i-1} q'_k\right| \quad \text{and}  \quad d^{(2)}=\max_{1\le i \le N}\left|\sum_{k=1}^{i} p'_k-\sum_{k=1}^{i-1} q'_k\right|,
%\end{eqnarray*}
%\begin{eqnarray*}
%d^{(3)}=\max_{1\le i \le n_1+n_2}\left|\sum_{k=1}^{i-1} p'_k-\sum_{k=1}^{i-1} q'_k\right| \quad \text{and}  \quad d^{(4)}=\max_{1\le i \le N}\left|\sum_{k=1}^{i} p'_k-\sum_{k=1}^{i-1} q'_k\right|,
%\end{eqnarray*}
%provided that $\sum_{k=1}^{0}q'_k=0$. Here $\theta_{(1)}< \ldots< \theta_{(N)}$ are the order statistics of $(\theta_k)_{1 \leq k \leq N}$,
%and $q'_k$ is the weight associated with $\theta_{(k)}$ such that for  each $1 \le k \le N$, $q'_k=q_j$ if $\theta_{(k)}=\theta_{j}$ for some $1 \le j \le N$.
%\EndLemma

 Next, we   show empirically that the proposed technique of goodness-of-fit test is robust against  prior specification the Dirichlet process' parameters (i.e., the concentration parameter and the base measure). To this end, we  have repeated Example 1 in Section 4 with two additional cases. In the first case,  we take $H$ to be the uniform distribution on $[0,1]$ and $a=1$. In the second case, we take  $H$ to be the standard normal distribution and $a=50$. The results are reported in Table 9. It follows clearly from Table 9 that the conclusions drawn in Example 1 are not affected by changing either $H$ or $a$. Specifically, fixing $a$ and changing $H$, will not change either the  distance or $U$. This is due to the fact that  the distance is measured  between two Dirichlet processes, which only depends on the cumulative weights of  each Dirichlet process. Thus, the distance and $U$ are independent from $H$. On the other hand, fixing $H$ and changing $a$ changes both the distance and $U$. If $a$ is increased (decreased), then the both distance and $U$ are decreased (increased) in such way that conclusions are not affected and remain the same.

\begin{table}[htbp]
  \centering
  \caption{Robustness of the Bayesian nonparametric test against changing the parameters of the Dirichlet process}
    \begin{tabular}{lccc}
    \hline
Samples &$d: H=N(0,1),a=1$  & $d: H=U[0,1],a=1$ & $d: H=N(0,1),a=50$  \\
    \hline
1&0.15 	&0.15   &  0.12 	    \\
2&0.42	&0.42	&  0.30		\\
3&0.26	&0.26	&  0.19 		\\
4&0.44	&0.44	&  0.30			\\
5&0.19	&0.19	&  0.15			\\
6&0.31	&0.31	&  0.21			\\
7&0.42	&0.42	&  0.29			\\
8&0.26	&0.26	&  0.19			\\
  \hline
$U$&0.20 &0.20      &  0.18       \\ \hline
    \end{tabular}%
  \label{tab:addlabel}%
\end{table}%

%\begin{table}[htbp]
%  \centering
%  \caption{Robustness of the Bayesian nonparametric test against changing the parameters of the Dirichlet process}
%    \begin{tabular}{lccccc}
%    \hline
% \multicolumn{1}{c}{\multirow{2}[4]{*}{Samples}} &\multicolumn{2}{c}{$H=U[0,1]$ and $a=1$}  & \multicolumn{1}{c}{} &\multicolumn{2}{c}{$H=U[0,1]$ and $a=1$} \\  \cline{2-3} \cline{5-6}   &$d$ & $U$  &    & $d$ & $U$ \\
%    \hline
%1	&0.16   &  &&0.12   &	    \\
%2	&0.81	&  &&0.55	&		\\
%3	&0.26	&  &&0.17   &		\\
%4	&0.35	&  &&0.24	&		\\
%5	&0.17	&  &&0.13	&		\\
%6	&0.28	&  &&0.20	&		\\
%7	&0.44	&  &&0.30	&		\\
%8	&0.27	&  &&0.19	&		\\
%  \hline
%    \end{tabular}%
%  \label{tab:addlabel}%
%\end{table}%

\section{Concluding Remarks}
A method  based on the Kolmogorov distance and approximate samples from the Dirichlet process is proposed to asses the equality of two unknown distributions. The new approach is simple and flexible such that it can be applied to any two-sample problem with any sample size. As shown in the power study, the proposed test has distinctly higher power than the power of several popular tests in many cases. In particular, the proposed test dominates the standard Kolmogorov-Smirnov test in all the cases examined in the present paper.   On the other hand, unlike most frequentist tests, the proposed test is not based on computing p-values. The main concern about using   p-values in testing statistical hypothesis is that they overestimate the evidence against the null hypothesis (Masson, 2011; Sawrtz, 1999; Wagenmakers, 2007). For instance, see Table 1.

The current study may lead to further research directions. For instance, it would be interesting to study the effect of selecting other distances such as the Wasserstein (or Kantorovich) distance and the Kullback-Leibler distance on the proposed approach. Another important extension is the generalization of the approach to construct a goodness-of-fit test for multivariate distributions. In principle, there is no need to change the methodology. However, the calculation of the distance  requires amendment in this case. Extending the approach to multivariate distributions will bypass the distribution-free problem for  the tests that  rely on the empirical distribution function. Finally,  similar to the frequentist's Kolmogorov-Smirnov test, it is possible to construct a test based on the fact that the two independent processes $\sqrt{m_1}\left(P^{*}_{m_1}-H^*_{m_1}\right)$ and $\sqrt{m_2}\left(P^{*}_{m_2}-H^{*}_{m_2}\right)$ converge jointly in distribution to the two independent Brownian bridges $B_{F}$ and $B_{G}$, where $F$ and $G$ are the ``true" distributions generating the data. Recall that, for a collection of Borel sets $\mathscr{S}$ in $\mathbb{R}$, a Gaussian process $\left\{B_F(S): S \in \mathscr{S}\right\}$ is called a {\emph{Brownian bridge with parameter measure }$F$  if ${E}\left[B_F(S)\right]=0$ for any $S \in \mathscr{S}$ and
${Cov}\left(B_F(S_1),B_F(S_2)\right)=F(S_1 \cap S_2)- F(S_1)F(S_2)$ for any  $S_1, S_2 \in \mathscr{S}$ (Kim and Bickel, 2003). We leave this direction for future work.

%\section{Acknowledgments} We thank the  Editor, Associate Editor and two anonymous referees  for their careful review of the paper and their extremely helpful comments.  Research of the third author is supported by the \textbf{Natural Sciences and Engineering Research Council of Canada (NSERC)}.

\begin{landscape}
 \section*{Appendix}
\begin{table}[h!]
  \centering
\renewcommand{\arraystretch}{1.1}
  \caption{Values of $U$ (upper bound of the 95\% prediction interval for the prior distance $d_0$).}
    \begin{tabular}{p{8pt}|	p{14pt}	p{14pt}	p{14pt}	p{14pt}	p{14pt}	p{14pt}	p{14pt}	p{14pt}	p{14pt}	p{14pt}	p{14pt}	p{14pt}	p{14pt}	p{14pt}	p{14pt}	p{14pt}	p{14pt}	p{14pt}	p{14pt}	p{14pt}
}
          & \multicolumn{20}{l}{$m_1$} \\
        $m_2$ & 1     & 2     & 3     & 4     & 5     & 6     & 7     & 8     & 9     & 10    & 11    & 12    & 13    & 14    & 15    & 16    & 17    & 18    & 19    & 20 \\ \hline
    1     & 0.93 &       &       &       &       &       &       &       &       &       &       &       &       &       &       &       &       &       &       &  \\
    2     & 0.90 & 0.86 &       &       &       &       &       &       &       &       &       &       &       &       &       &       &       &       &       &  \\
    3     & 0.87 & 0.82 & 0.79 &       &       &       &       &       &       &       &       &       &       &       &       &       &       &       &       &  \\
    4     & 0.87 & 0.80 & 0.78 & 0.74 &       &       &       &       &       &       &       &       &       &       &       &       &       &       &       &  \\
    5     & 0.84 & 0.79 & 0.74 & 0.72 & 0.69 &       &       &       &       &       &       &       &       &       &       &       &       &       &       &  \\
    6     & 0.84 & 0.78 & 0.74 & 0.71 & 0.67 & 0.66 &       &       &       &       &       &       &       &       &       &       &       &       &       &  \\
    7     & 0.82 & 0.78 & 0.72 & 0.69 & 0.67 & 0.64 & 0.63 &       &       &       &       &       &       &       &       &       &       &       &       &  \\
    8     & 0.82 & 0.75 & 0.69 & 0.68 & 0.64 & 0.64 & 0.62 & 0.60 &       &       &       &       &       &       &       &       &       &       &       &  \\
    9     & 0.79 & 0.76 & 0.71 & 0.69 & 0.64 & 0.61 & 0.61 & 0.60 & 0.58 &       &       &       &       &       &       &       &       &       &       &  \\
    10    & 0.81 & 0.76 & 0.71 & 0.68 & 0.62 & 0.60 & 0.58 & 0.58 & 0.58 & 0.55 &       &       &       &       &       &       &       &       &       &  \\
    11    & 0.83 & 0.74 & 0.67 & 0.63 & 0.61 & 0.59 & 0.58 & 0.57 & 0.55 & 0.54 & 0.52 &       &       &       &       &       &       &       &       &  \\
    12    & 0.79 & 0.73 & 0.66 & 0.67 & 0.63 & 0.58 & 0.58 & 0.56 & 0.55 & 0.54 & 0.52 & 0.52 &       &       &       &       &       &       &       &  \\
    13    & 0.80 & 0.72 & 0.67 & 0.64 & 0.59 & 0.58 & 0.55 & 0.56 & 0.53 & 0.54 & 0.51 & 0.51 & 0.50 &       &       &       &       &       &       &  \\
    14    & 0.79 & 0.73 & 0.66 & 0.64 & 0.60 & 0.57 & 0.57 & 0.54 & 0.53 & 0.52 & 0.51 & 0.51 & 0.50 & 0.49 &       &       &       &       &       &  \\
    15    & 0.79 & 0.71 & 0.67 & 0.62 & 0.60 & 0.57 & 0.55 & 0.54 & 0.51 & 0.51 & 0.50 & 0.50 & 0.50 & 0.47 & 0.47 &       &       &       &       &  \\
    16    & 0.79 & 0.70 & 0.65 & 0.61 & 0.59 & 0.59 & 0.56 & 0.53 & 0.53 & 0.52 & 0.51 & 0.51 & 0.49 & 0.47 & 0.46 & 0.45 &       &       &       &  \\
    17    & 0.79 & 0.72 & 0.67 & 0.63 & 0.59 & 0.56 & 0.56 & 0.53 & 0.52 & 0.50 & 0.49 & 0.49 & 0.47 & 0.47 & 0.46 & 0.45 & 0.45 &       &       &  \\
    18    & 0.80 & 0.72 & 0.66 & 0.62 & 0.60 & 0.58 & 0.55 & 0.53 & 0.52 & 0.51 & 0.48 & 0.47 & 0.47 & 0.48 & 0.45 & 0.45 & 0.44 & 0.44 &       &  \\
    19    & 0.79 & 0.70 & 0.65 & 0.63 & 0.58 & 0.57 & 0.53 & 0.54 & 0.51 & 0.49 & 0.49 & 0.47 & 0.46 & 0.46 & 0.44 & 0.45 & 0.45 & 0.42 & 0.44 &  \\
    20    & 0.79 & 0.73 & 0.65 & 0.62 & 0.59 & 0.57 & 0.54 & 0.51 & 0.49 & 0.49 & 0.49 & 0.47 & 0.46 & 0.44 & 0.44 & 0.45 & 0.43 & 0.43 & 0.42 & 0.42 \\ \hline
    \end{tabular}%
  \label{tab:addlabel}%
\end{table}%
\end{landscape}
\restoregeometry

%%%%%%%%%%%%%%%%%%%%%%%%%%%%%%%%%%%%%%%%%%%%%%%%%%%%%%%%%%%%%%%%%%%%%%%%%%%%%%%%%%%%%%%%%%%%%%%%%%%%%%%%%%%%%%%%%%%%%%%%%%%%
\end{document}